\RequirePackage{ifpdf}
\ifpdf 
\documentclass[pdftex]{sigma}
\else
\documentclass{sigma}
\fi

\def\Mg{{\overline{M}_{g,n}(\mathbb{P}^{r},d)}}
\def\M0{{\overline{M}_{0,n}(\mathbb{P}^{r},d)}}
\def\mm{{\overline{\mathcal{M}}}}
\def\m0{{\overline{\mathcal{M}}_{0,n}(\mathbb{P}^{r},d)}}
\def\mznXb{{\overline{\mathcal{M}}_{0,n}(X,\beta)}}

\def\mg{{\overline{\mathcal{M}}_{g,n}(\mathbb{P}^{r},d)}}
\def\nobarm0{{\mathcal{M}_{0,n}(\mathbb{P}^{r},d)}}

\def\gr{\operatorname{gr}}

\def\ser{\operatorname{Serre}}
\def\PGL{\operatorname{PGL}}
\def\SL{\operatorname{SL}}
\def\Aut{\operatorname{Aut}}
\def\Trace{\operatorname{Tr}}

\def\one{\textrm{\makebox[0.02in][l]{1}1}}

\begin{document}

\allowdisplaybreaks

\renewcommand{\textfraction}{0.01}
\renewcommand{\topfraction}{0.99}

\newcommand{\chews}[2]{\genfrac{[}{]}{0pt}{}{#1}{#2}}
\newcommand{\lchoose}[2]{\dbinom{#1}{#2}}

\renewcommand{\PaperNumber}{085}

\FirstPageHeading

\ShortArticleName{An Additive Basis for the Chow Ring of $\mm_{0,2}({\mathbb P}^r,2)$}

\ArticleName{An Additive Basis for the Chow Ring of $\boldsymbol{\mm_{0,2}({\mathbb P}^r,2)}$}

\Author{Jonathan A. COX}

\AuthorNameForHeading{J.A. Cox}

\Address{Department of Mathematical Sciences, SUNY Fredonia, Fredonia, New York 14063, USA}

\Email{\href{mailto:Jonathan.Cox@fredonia.edu}{Jonathan.Cox@fredonia.edu}}

\URLaddress{\url{http://www.fredonia.edu/faculty/math/JonathanCox/}}

\ArticleDates{Received July 03, 2007, in f\/inal form August 28, 2007; Published online August 31, 2007}

\Abstract{We begin a study of the intersection theory of the moduli spaces of degree two stable maps from
two-pointed rational curves to arbitrary-dimensional projective space. First we compute the
Betti numbers of these spaces using Serre polynomial and equivariant Serre polynomial methods
developed by E.~Getzler and R.~Pandharipande. Then, via the excision sequence, we compute an
additive basis for their Chow rings in terms of Chow rings of nonlinear Grassmannians, which
have been described by Pandharipande. The ring structure of one of these Chow rings is addressed
in a sequel to this paper.
}


\Classification{14C15; 14D22}

\section{Introduction}
\label{sec:intro}

Let $\mg$ be the moduli space of stable maps from $n$-pointed,
genus $g$ curves to ${\mathbb P}^r$ of degree $d$.
In this article we will begin a study of the intersection theory of the moduli
spaces $\mm_{0,2}({\mathbb P}^r,2)$.

Moduli spaces of stable maps have proven useful in studying
both superstring theory and enumerative geometry.
Developing a solid mathematical foundation for computation of certain
numbers in string theory
was the primary motivation behind the introduction of moduli spaces of
stable maps in \cite{KM}. Examples include the {\em instanton numbers},
which intuitively count the number of
holomorphic instantons on a space $X$ (nonconstant holomorphic maps from Riemann
surfaces to $X$). Instanton numbers are calculated using other values
called {\em Gromov--Witten invariants}.
Naively, Gromov--Witten invariants should count the number of curves of
a certain homology class and genus which pass through certain subvarieties of
the target space. More specif\/ically,
let $X$ be a projective manifold, $\beta\in H_2(X)$, $g$ and $n$ nonnegative
integers. Let $\gamma_1,\ldots,\gamma_n\in H^*(X)$ be cohomology classes such that
there exist subvarieties $Z_1,\ldots,Z_n$
with $Z_i$ representing the Poincar\'{e} dual of $\gamma_i$.
Then the Gromov--Witten invariant $\langle\gamma_1,\ldots,\gamma_n\rangle_{g,\beta}$ should
count
the number of genus $g$ curves of class $\beta$ that intersect all of the $Z_i$.
Also of interest are {\em gravitational correlators},
which generalize Gromov--Witten invariants. Gravitational correlators are def\/ined and computed mathematically as
intersection numbers on the moduli space of stable maps.
See \cite[Chapter~10]{CK} for a rigorous def\/inition and the sequel \cite{C} for some computations that f\/low from
the results of this paper.

In the dozen years since Kontsevich introduced the concept in \cite{KM} and
\cite{Ko}, the moduli space of stable maps has been exploited to solve
a plethora of enumerative problems for curves.
As a rule, these results were derived without a
complete description of the Chow rings involved. Instead, the requisite
intersection numbers were calculated somewhat indirectly,
most often using the method of {\em localization}. (An overview of localization
is given
in \cite[Chapter 9]{CK}.)  Such a complete description is the key step
in giving another, more direct, computation of these enumerative numbers
and possibly many others.
Since a presentation for a ring gives
an easy way to compute all products in the ring, giving
presentations for the Chow rings of moduli spaces of stable maps
is the clear path toward attaining a full and direct knowledge of their
intersection theory.
As a consequence, this also helps give a new and more direct way of
determining values of instanton numbers, Gromov--Witten invariants, and
gravitational correlators.

Until recently, presentations for Chow rings of moduli spaces of stable maps were
known only in a few special cases. Most of these had projective space as
the target of the stable maps, and in this case the moduli space
$\mg$ depends
on four nonnegative integer parameters: the genus $g$ of the curves, the
number $n$ of marked points on the curves, the dimension $r$ of the target
projective space, and the degree $d$ of the stable maps.
Most impressive was the presentation of
$A^*(\mm_{0,1}({\mathbb P}^r,d))$ for arbitrary $d$ and $r$ described by Mustata and Mustata in \cite{MM}.
  Behrend and
O'Halloran gave a presentation for $A^*(\mm_{0,0}({\mathbb P}^r,2))$ and
conjectured a presentation for $A^*(\mm_{0,0}({\mathbb P}^r,3))$ in
\cite{BO}. Also of relevance, Oprea more recently described a system
of {\em tautological subrings} of the cohomology (and hence Chow)
rings in the genus zero case and showed that, if the target $X$ is
an $\SL$ f\/lag variety, then all rational cohomology classes on
$\mznXb$ are tautological. This gave, at least in principle, a set
of generators for any such Chow ring, namely its tautological
classes. He furthermore described an additive generating set for the
cohomology ring of any genus zero moduli space (with target a
projective algebraic variety). Finally, he speculated that all
relations between the tautological generators are consequences of
the topological recursion relations. These developments were
substantial steps toward describing presentations for the Chow rings
of moduli spaces of stable maps in much more general cases. See
\cite{O} and~\cite{O2} for more details. More basic examples include
$A^*(\mm_{0,n}({\mathbb P}^r,0))\simeq A^*({\mathbb P}^r)\times A^*({\overline{M}}_{0,n})$, where
${\overline{M}}_{0,n}$ is the moduli space of stable curves.  This case reduces
to f\/inding presentations for the rings $A^*({\overline{M}}_{0,n})$, and Keel
did so in \cite{K}. Also, $\mm_{0,0}({\mathbb P}^r,1)$ is isomorphic to
${\mathbb G}(1,r)$, the Grassmannian of lines in projective space, and
$\mm_{0,1}({\mathbb P}^r,1)$ is isomorphic to ${\mathbb F}(0,1;r)$, the f\/lag variety
of pointed lines in projective space. The
spaces $\mm_{0,n}({\mathbb P}^1,1)$ are Fulton-MacPherson compactif\/ications
of conf\/iguration spaces of ${\mathbb P}^1$. Presentations for their Chow
rings were given by Fulton and MacPherson in \cite{FM}. Detailed
descriptions of Chow rings of spaces $\mg$, with $g>0$, are almost
nonexistent (although some progress is now being made for $g=1$). Additional complications arise in this case.

This was the state of af\/fairs up to the posting of this article, which
lays the foundation for computing presentations of
the Chow rings of the spaces
$\mm_{0,2}({\mathbb P}^r,2)$. This computation is completed for the case $r=1$
in \cite{C}, the sequel to this article, where we obtain the presentation
\[
A^*(\mm_{0,2}({\mathbb P}^1,2))\simeq\frac{{\mathbb Q}[D_0,D_1,D_2,H_1,H_2,\psi_1,\psi_2]}
{\left(\!\begin{array}{c} H_1^2, H_2^2,D_0\psi_1,D_0\psi_2,D_2-\psi_1-\psi_2,
\psi_1-\frac{1}{4}D_1-\frac{1}{4}D_2-D_0+H_1, \vspace{1mm}\\
\psi_2-\frac{1}{4}D_1-\frac{1}{4}D_2-D_0+H_2, (D_1+D_2)^3,
D_1\psi_1\psi_2 \end{array}\!\right)}.
\]
This gave the f\/irst known presentation for a Chow ring of a
moduli space of stable maps of degree greater than one with more than one
marked point.
The sequel also employs the presentation to give a new
computation of the genus zero, degree two, two-pointed gravitational
correlators of ${\mathbb P}^1$. Algorithms for computing theses values have previously
been developed; see \cite{KM2} and \cite{CK}, for example.

Knowing the Betti numbers of $\mm_{0,2}({\mathbb P}^r,2)$ is the
important f\/irst step in our computation, since it will give us a good idea of
how many generators and relations to expect in each degree.  We accomplish
this in Section \ref{sec:ser} by using the equivariant Serre polynomial method
of Getzler and Pandharipande. I owe much gratitude to Getzler and Pandharipande
for supplying a copy of
their unpublished preprint \cite{GP}, which
provided the inspiration for Section \ref{sec:ser}.
In Section \ref{sec:basis}, we use the presentations of
the Chow rings for the moduli spaces ${\mathcal{M}}_{0,0}({\mathbb P}^r,d)$ from \cite{P4}
together with excision to determine a generating set for the Chow ring
$A^*(\mm_{0,2}({\mathbb P}^r,2))$. Comparing the dimension of each graded piece
with the Betti numbers from Section \ref{sec:ser}, we conclude that this generating
set is actually an additive basis.

Since the completion of this article, additional results have been announced.
Getzler and Pandharipande have substantially reworked
the material from and completed the computations begun in their preprint, giving Betti numbers
for all the spaces $\m0$ in \cite{GP2}. However, they do this via an indirect method, using a
generalization of the Legendre transform, and the result is an algorithm for
computing the Betti numbers rather than a closed formula for them. The approach taken here is
more straightforward and sheds more light on the geometry involved. Finally, Mustata and Mustata extended
their results for one-pointed spaces to supply
presentations for the Chow rings of all spaces $\m0$ in \cite{MM2}. Again, the methods of the present paper are
simpler, and the presentation obtained in the sequel is more explicit.

\subsection{Conventions and preliminary comments}
\label{sec:mod}

We will work over the f\/ield ${\mathbb C}$ of complex numbers.
Let $\underline{n}={\mathbb N}\cap[1,n]$ be the
initial segment consisting of the f\/irst $n$ natural numbers.

The moduli stack $\mg$
captures all the data of the moduli problem for stable maps, while the moduli scheme
$\Mg$ loses some information, including that of
automorphisms of families. Since retaining
all of this data leads to a more beautiful, powerful, and complete theory,
we will work with the stack incarnations of the moduli spaces rather than the
coarse moduli schemes.

Let $H^*(F)$ denote the rational de Rham cohomology ring of a Deligne--Mumford
stack $F$.

\begin{proposition}[Homology isomorphism]
\label{hi}
Let $X$ be a flag variety. Then
there is a canonical ring isomorphism
\begin{gather}\label{a2h}
A^*(\mznXb)\rightarrow H^*(\mznXb).
\end{gather}
\end{proposition}

 See \cite{O} for a proof. Following
\cite{K}, we call any scheme or stack $Y$ an {\em HI scheme} or {\em
HI stack} if the canonical map $A^*(Y)\rightarrow H^*(Y)$ is an isomorphism.
In particular, by Proposition 1 $Y=\m0$ is an HI stack since ${\mathbb P}^r$ is a f\/lag
variety. This allows us to switch freely between cohomology and Chow
rings. We should note that the isomorphism doubles degrees. The
degree of a $k$-cycle in $A^k(\mznXb)$, called the algebraic degree,
is half the degree of its image in $H^{2k}(\mznXb)$.

\section[The Betti numbers of $\mm_{0,2}({\mathbb P}^r,2)$]{The Betti numbers of $\boldsymbol{\mm_{0,2}({\mathbb P}^r,2)}$}
\label{sec:ser}

\subsection[Serre polynomials and the Poincar\'{e} polynomial
of $\mm_{0,2}({\mathbb P}^r,2)$]{Serre polynomials and the Poincar\'{e} polynomial
of $\boldsymbol{\mm_{0,2}({\mathbb P}^r,2)}$}
\label{sec:poincare}

This section owes much to Getzler and Pandharipande, who provide the
framework for computing the Betti numbers of all the spaces $\m0$
in \cite{GP}. However, we will take the def\/initions and basic results from
other sources, and prove their theorem in the special case that we need.
We will compute a formula for the Poincar\'{e} polynomials of the moduli
spaces $\mm_{0,2}({\mathbb P}^r,2)$ using what are called Serre polynomials in
\cite{GP} and Serre characteristics in \cite{GP2}. (These polynomials are also known as virtual
Poincar\'{e} polynomials or E-polynomials.)
Serre polynomials are def\/ined for
varieties over ${\mathbb C}$ via the mixed Hodge theory of Deligne
(\cite{D}). Serre conjectured the existence of polynomials satisfying
the key properties given below. A formula was later given by
Danilov and Khovanski\u{\i} in \cite{DK}. If $(V,F,W)$ is a mixed Hodge
structure over ${\mathbb C}$, set
\[V^{p,q}=F^p\gr_{p+q}^WV\cap \bar{F}^q\gr_{p+q}^WV
\]
and let ${\mathcal X}(V)$ be the Euler characteristic of $V$ as a graded vector space.
Then
\[\ser(X)=\sum_{p,q=0}^\infty u^pv^q{\mathcal X}(H_c^\bullet(X,{\mathbb C})^{p,q}).
\]
If $X$ is a smooth projective variety, then the
Serre polynomial of $X$ is just its Hodge polynomial:
\[\ser(X)=\sum_{p,q=0}^\infty(-u)^p(-v)^q\dim H^{p,q}(X,{\mathbb C}).\]
If $X$ further satisf\/ies $H^{p,q}(X,{\mathbb C})=0$ for $p\neq q$, then we can
substitute a new variable $q=uv$ for $u$ and $v$.
In this case, the coef\/f\/icients of the Serre polynomial of $X$ give its
Betti numbers, so that $\ser(X)$ is the Poincar\'{e} polynomial of $X$.

We will use two additional key properties of Serre
polynomials. The f\/irst gives a compatibility
with decomposition: If $Z$ is a closed subvariety of $X$, then
$\ser(X)=\ser(X\backslash Z)+\ser(Z)$.  Second, it respects products:
$\ser(X\times Y)=\ser(X)\ser(Y)$.  (The latter is actually a~consequence of the previous
properties.) It follows from these two properties
that the Serre polynomial
of a f\/iber space is the product of the Serre polynomials of the base and the
f\/iber. The def\/inition and properties above come from \cite{Ge}.
We also use the following consequence of the Eilenberg--Moore spectral
sequence, which is essentially Corollary 4.4 in~\cite{Sm}.

\begin{proposition}
Let $Y\rightarrow B$ be a fiber space with $B$ simply connected, and let $X\rightarrow B$
be continuous. If $H^*(Y)$ is a free $H^*(B)$-module, then
\[H^*(X\times_B Y)\simeq H^*(X)\otimes_{H^*(B)}H^*(Y)\]
as an algebra.
\end{proposition}

Since we deal exclusively with cases where the isomorphism (\ref{a2h})
holds, there is never any torsion in the cohomology. Thus we have the
following.

\begin{corollary}
\label{serfiber}
Let $X$ and $Y$ be varieties over a simply connected base $B$, and suppose
either $X$ or $Y$ is locally trivial over $B$. Then
\[\ser(X\times_B Y) = \frac{\ser(X)\ser(Y)}{\ser(B)}.\]
\end{corollary}

We will sometimes use the notation $Y/B$ for the f\/iber of a f\/iber space
$Y\rightarrow B$.

To extend this setup to Deligne--Mumford stacks, where automorphism
groups can be nontrivial (but still f\/inite),
{\em equivariant} Serre polynomials
are needed.  Let $G$ be a f\/inite group acting on a 
variety $X$.
The idea is this: The action  of $G$ on $X$
induces an action on its cohomology (preserving the mixed Hodge
structure), which in turn gives a representation of $G$ on each (bi)graded
piece of the cohomology. The cohomology of the quotient variety $X/G$,
and hence of
the quotient stack $[X/G]$, is the part of the cohomology of $X$ which
is f\/ixed by the $G$-action, {\em i.e.}, in each degree the subspace on which
the representation is trivial.

Our def\/inition comes from
\cite{Ge}.
The equivariant Serre polynomial $\ser(X,G)$ of $X$ is given by the formula
\[\ser_g(X)=\sum_{p,q=0}^\infty u^pv^q
\sum_i(-1)^i\Trace(g|(H_c^i(X,{\mathbb C}))^{p,q}).
\]
for each element $g\in G$.
We can also describe the equivariant Serre polynomial more compactly
with the formula
\[\ser(X,G)=\sum_{p,q=0}^\infty u^pv^q
\sum_i(-1)^i[H_c^i(X,{\mathbb C})^{p,q}],\]
taken from \cite{GP}.
In the
case $G=S_n$, we write $\ser_n(X)$ for $\ser(X,S_n)$.
A $G$-equivariant Serre polynomial takes values
in $R(G)[u,v]$, where $R(G)$ is the virtual representation ring of $G$.
The augmentation morphism $\epsilon:R(G)\rightarrow{\mathbb Z}$, which extracts the
coef\/f\/icient of the trivial representation~$\one$ from an element
of $R(G)$, extends to an
augmentation morphism $R(G)[u,v]\rightarrow{\mathbb Z}[u,v]$. If $G$ acts on
a quasi-projective variety $X$, the Serre polynomial of the
quotient stack $[X/G]$ is the
augmentation of the equivariant Serre polynomial of $X$.

Every virtual representation ring $R(G)$
has the extra structure
of a {\em $\lambda$-ring}. See \cite{Knutson} for the def\/initions and basic
properties of $\lambda$-rings and pre-$\lambda$-rings. Here we just brief\/ly state
the most relevant facts.

Let $V$ be a $G$-module. Then $\lambda_i(V)$ is the {\em $i$'th exterior power}
$\Lambda^iV$
of $V$, where we def\/ine $g\in G$ to act by
$g(v_1\wedge\cdots\wedge v_i)=gv_1\wedge\cdots\wedge gv_i$. Def\/ine $\lambda_0(V)$
to be the trivial one-dimensional representation.
(One can similarly def\/ine a $G$-module structure
on the $i$'th symmetric power~$S^iV$.) Knutson proves in
\cite[Chapter II]{Knutson} that
these exterior power operations give $R(G)$ the structure of a $\lambda$-ring
for any f\/inite group $G$. Addition is given by $[V]+[W]=[V\oplus W]$, and the
product is $[V]\cdot[W]=[V\otimes W]$, both with the naturally induced actions.

Knutson also shows that ${\mathbb Z}$ is a $\lambda$-ring with $\lambda$-operations
given via $\lambda_t(m)=(1+t)^m$, where by def\/inition $\lambda_t(m)=\sum \lambda_i(m)t^i$.
For $m,n\geq 0$, this gives
$\lambda_n(m)=\lchoose{m}{n}$. Finally, he shows that
if $R$ is a $\lambda$-ring, then there is
a unique structure of $\lambda$-ring on $R[x]$ under which
$\lambda_k(rX^n)=\lambda_k(r)X^{nk}$ for $n,k\in{\mathbb N}\cup\{0\}$ and $r\in R$.
This gives a $\lambda$-ring structure on ${\mathbb Z}[q]$. The augmentation morphism
$\epsilon:R(G)\rightarrow{\mathbb Z}$ is a {\em map of $\lambda$-rings}; it commutes with the
$\lambda$-operations.

We will use the following facts about Serre polynomials and equivariant
Serre polynomials.
For $n\in{\mathbb N}$, let $[n]=\frac{q^n-1}{q-1}$. Then $[n+1]$ is the Serre
polynomial of ${\mathbb P}^n$, as is clear from the presentation for its Chow ring.
Getzler and Pandharipande prove that the
Serre polynomial of the Grassmannian $G(k,n)$ of
$k$-planes in ${\mathbb C}^n$ is the $q$-binomial coef\/f\/icient
\[\chews{n}{k}=\frac{[n]!}{[k]![n-k]!},\]
where $[n]!=[n][n-1]\cdot\cdot\cdot[2][1]$. We
will prove this formula in the special case $k=2$.

\begin{lemma}\label{grasser}
The Serre polynomial of $G(2,n)$ is $\chews{n}{2}$.
\end{lemma}

\begin{proof}
We can work with the Grassmannian
${\mathbb G}(1,n-1)$ of lines in ${\mathbb P}^{n-1}$ since $G(2,n)\simeq{\mathbb G}(1,n-1)$.
The universal ${\mathbb P}^1$-bundle over ${\mathbb G}(1,n-1)$ is isomorphic to
${\mathbb F}(0,1;n-1)$, the f\/lag variety of pairs $(p,\ell)$ of a point $p$
and a line $\ell$ in ${\mathbb P}^{n-1}$ with $p\in\ell$. On the other hand,
there is a~projection ${\mathbb F}(0,1;n-1)\rightarrow{\mathbb P}^{n-1}$ taking $(p,\ell)$
to $p$. Its f\/iber over a point $p$ is $\{\ell\, |\, p\in\ell\}$,
which is isomorphic to ${\mathbb P}^{n-2}$. (To see this isomorphism, f\/ix a
hyperplane $H\subset{\mathbb P}^{n-1}$ not containing $p$ and map each line to
its intersection with $H$.) It follows that $\ser({\mathbb F}(0,1;n-1))
=[n][n-1]$. Since $\ser({\mathbb F}(0,1;n-1))=\ser({\mathbb G}(1,n-1))[2]$ also, we
are able to conclude that $\ser({\mathbb G}(1,n-1))=[n][n-1]/[2]$.
\end{proof}

Next, since $\PGL(2)$
is the complement of a quadric surface in ${\mathbb P}^3$,
$\ser(\PGL(2))=[4]-[2]^2=q^3-q$.

In addition to the $\lambda$-operations, every $\lambda$-ring $R$
has {\em $\sigma$-operations}
as well. These can be def\/ined
in terms of the $\lambda$-operations by $\sigma_k(x)=(-1)^k\lambda_k(-x)$.
Routine checking shows that the $\sigma$-operations also give $R$ the structure of
a pre-$\lambda$-ring.
Here
we simply note the following formulas for the $\lambda$-ring ${\mathbb Z}[q]$
\[\sigma_k([n])=\chews{n+k-1}{k} \qquad \text{and}\qquad
\lambda_k([n])=q^{k \choose 2}\chews{n}{k}.\]
Proofs of these formulas can be found in
\cite[Section I.2]{Mac}. Next we explain why these formulas are
relevant. Let $\epsilon$ be the sign representation of $S_n$. Note the
identity $\epsilon^2=\one$. We will prove the following claim from
\cite{GP}.

\begin{lemma}\label{ser2X2}
If $X$ is a smooth variety and $S_2$ acts on $X^2$ by switching the factors, then
\[\ser_2(X^2)=\sigma_2(\ser(X))\text{\em \one}+\lambda_2(\ser(X))\epsilon.\]
\end{lemma}

\begin{proof}
Let $V$ be a vector space. Now $V\otimes V=S^2V\oplus \Lambda^2V$ as $S_2$-modules,
with $S_2$ acting by switching the factors of $V\otimes V$,
trivially on $S^2V$, and by sign on $\Lambda^2V$.
If 0 is the zero representation, certainly $\lambda_i(0)=0$ for $i>0$.
We use this fact and the properties of $\lambda$-rings
to obtain
\begin{gather*}
0  =  \lambda_2([V]-[V])
 =  \one\cdot\lambda_2(-[V])+[V]\cdot(-[V])+\lambda_2[V]\cdot\one\\
\phantom{0}   =  \lambda_2(-[V])-[S^2V]-[\Lambda^2V]+\lambda_2[V].
\end{gather*}

Since $\sigma_2[V]=\lambda_2(-[V])$, this implies $\sigma_2[V]=[S^2V]$. Since $X$ is
smooth, $H^*(X^2)=H^*(X)\otimes H^*(X)$, with the action of $S_2$ switching the
factors. Applying the above with $V=H^*(X)$ gives
$[H^*(X^2)]=\sigma_2[H^*(X)]+\lambda_2[H^*(X)]$. Breaking this down by
(cohomological) degree,
we have $[H^i(X^2)]q^i=[\sigma_2[H^*(X)]]_iq^i+[\lambda_2[H^*(X)]]_iq^i$.
We need to show
that $[H^i(X^2)]q^i=[\sigma_2(\ser(X))]_i\one +[\lambda_2(\ser(X))]_i\epsilon$. We will show
the equality of
the f\/irst summands of each expression; showing equality of
the terms involving $\lambda_2$ is easier. First, by induction the identity
$\lambda_2(-m)=\lambda_2(m+1)$ holds. Second, note that any pre-$\lambda$-operation $\lambda_2$
acts on sums by $\lambda_2(\sum_i x_i)=\sum_i \lambda_2(x_i)+\sum_{i<j}x_ix_j$.
Third, note that vector spaces in the following computation live in the
graded algebra $H^*(X)\otimes H^*(X)$, and we will apply the usual rules for
grading in a~tensor product. Finally, all of the representations below
are trivial.
We f\/ind
\begin{gather*}
 [\sigma_2[H^*(X)]]_iq^i  =  \left[\sigma_2\left[\sum_j H^j(X)\right]\right]_iq^i =  \left[\sum_j[S^2H^j(X)]+\sum_{j<k}[H^j(X)\otimes H^k(X)]\right]_iq^i\\
\phantom{[\sigma_2[H^*(X)]]_iq^i}{} = \begin{cases}
\Bigg([S^2H^{i/2}(X)]+\sum\limits_{\stackrel{\scriptstyle j+k=i}{j<k}}
[H^j(X)\otimes H^k(X)]\Bigg)q^i
& \text{if $i$ is even,}\\
\Bigg(\sum\limits_{\stackrel{\scriptstyle j+k=i}{j<k}}[H^j(X)\otimes H^k(X)]\Bigg)q^i
& \text{if $i$ is odd,}
\end{cases}\\
\phantom{[\sigma_2[H^*(X)]]_iq^i}{} = \begin{cases}
\Bigg(\lchoose{{\displaystyle h^{i/2}(X)+1}}{{\displaystyle 2}}\one
+\sum\limits_{\stackrel{\scriptstyle j+k=i}{j<k}}h^j(X)h^k(X)\one\Bigg)q^i
& \text{if $i$ is even,}\\
\Bigg(\sum\limits_{\stackrel{\scriptstyle j+k=i}{j<k}}h^j(X)h^k(X)\one\Bigg)q^i
& \text{if $i$ is odd.}
\end{cases}
\end{gather*}
On the other hand,
\begin{gather*}
  [\sigma_2(\ser(X))]_i\one
=  [\lambda_2(-\sum h^j(X)q^j)]_i\one \\
\phantom{[\sigma_2(\ser(X))]_i\one}{} =  \left[\sum \lambda_2(-h^j(X))q^{2j}+\sum_{j<k}h^j(X)h^k(X)q^{j+k}
    \right]_i\one \\
\phantom{[\sigma_2(\ser(X))]_i\one}{} = \begin{cases}
\displaystyle\Bigg(\lchoose{{\displaystyle h^{i/2}(X)+1}}{{\displaystyle 2}}q^i
+\sum\limits_{\stackrel{\scriptstyle j+k=i}{j<k}}h^j(X)h^k(X)q^i\Bigg)\one
& \text{if $i$ is even,}\\
\displaystyle \Bigg(\sum\limits_{\stackrel{\scriptstyle j+k=i}{j<k}}h^j(X)h^k(X)q^i\Bigg)\one
& \text{if $i$ is odd.}
\end{cases}  \tag*{\raisebox{-1cm}{\qed}}
\end{gather*}
\renewcommand{\qed}{}
\end{proof}

As a corollary, the ordinary Serre polynomial of $[X^2/S_2]$ is
$\sigma_2(\ser(X))$.

The following proposition gives a key fact used in our computations. Notice that it
refers to the locus ${\mathcal{M}}_{0,0}({\mathbb P}^r,d))$ of stable maps with
smooth domain curve, which is a proper (dense) subset of the compactif\/ied
moduli space $\mm_{0,0}({\mathbb P}^r,d)$.

\begin{proposition}
\label{Serre00}
If $d>0$, $\ser({\mathcal{M}}_{0,0}({\mathbb P}^r,d))=q^{(d-1)(r+1)}\chews{r+1}{2}$.
\end{proposition}

This follows from Pandharipande's proof in \cite{P4} that the
Chow ring of the nonlinear Grassmannian $M_{{\mathbb P}^k}({\mathbb P}^r,d)$ is
isomorphic to the Chow ring of the ordinary Grassmannian ${\mathbb G}(k,r)$.
If $k=1$, the nonlinear Grassmannian is ${\mathcal{M}}_{0,0}({\mathbb P}^r,d)$. (The Serre
polynomial grades by dimension rather than codimension. This is why the
shifting factor $q^{(d-1)(r+1)}$ appears.)

Recall that
${\mathcal{M}}_{0,n}({\mathbb P}^r,0)\simeq M_{0,n}\times{\mathbb P}^r$, so that the Serre polynomials
of these spaces are easy to compute.

Finally,
let $F(X,n)$ be the conf\/iguration space of $n$ distinct labeled points in
a nonsingular variety $X$. Fulton and MacPherson show in \cite{FM} that
\[\ser(F(X,n))=\prod_{i=0}^{n-1}(\ser(X)-i)
.\]

In order to compute
the Serre polynomial of a moduli space of stable maps,  we can stratify
it according to the degeneration types of the maps and compute the
Serre polynomial of each stratum separately.  The degeneration types
of maps in $\m0$
are in 1--1 correspondence with stable $(n,d)$-trees via taking the
dual graph of a stable map.  These concepts were def\/ined and developed by Behrend
and Manin in \cite{BM}.

We are now ready to compute the Poincar\'{e} polynomials of some moduli spaces
of stable maps.

\begin{proposition}
The Poincar\'{e} polynomial of $\mm_{0,2}({\mathbb P}^r,2)$ is
\begin{gather}\label{ser2r2}
\ser(\mm_{0,2}({\mathbb P}^r,2))=
\left(\sum_{i=0}^r q^i\right)\left(\sum_{i=0}^{r-1} q^i\right)
\left(\sum_{i=0}^{r+2} q^i+2\sum_{i=1}^{r+1} q^i+2\sum_{i=2}^{r} q^i\right),
\end{gather}
and the Euler characteristic of $\mm_{0,2}({\mathbb P}^r,2)$ is $r(r+1)(5r+3)$.
\end{proposition}

\begin{proof}
We begin by stratifying $\mm_{0,2}({\mathbb P}^r,2)$ according to the degeneration
type of the stable maps. Since the strata are locally closed,
the compatibility of Serre polynomials
with decomposition allows us to compute the Serre polynomial of each
stratum separately and add up the results to obtain
$\ser(\mm_{0,2}({\mathbb P}^r,2))$.

\looseness=-1
Each stratum is isomorphic to a f\/inite group quotient of a
f\/iber product of moduli spaces of stable maps from smooth domain curves via the following procedure.
Given a stable map $(C,x_1,x_2,f)$, consider
the normalization of $C$.
It consists of a disjoint union of smooth curves~$C_i$ corresponding to the
components of $C$, and there are maps $f_i$
from each curve to ${\mathbb P}^r$ naturally
induced by $f$.  Furthermore, auxiliary marked points are added to retain
data about the node locations. The result is a collection of stable maps with
smooth domain curves,
one for each component of $C$. The evaluations of auxiliary marked points
corresponding to the same node must agree. This gives rise to
a f\/iber product of
moduli spaces type $\nobarm0$, together with a~morphism onto the stratum coming from the normalization map.
There can be automorphisms of the stable maps in the stratum that are not
accounted for by the f\/iber product. These occur when there is a collection
of connected unions $U_i$ of components that satisfy the following conditions:
\begin{enumerate}

\item[(1)] none of the $U_i$ contain any marked points;

\item[(2)] restrictions of $f$ to $U_i$ and $U_j$ give isomorphic maps
for all $i$ and $j$.
\end{enumerate}

 These automorphisms correspond exactly to the
automorphisms of the dual graph $\Gamma$ of the stratum. Proving the assertion that
$\Aut(\Gamma)$ is the right group to quotient by appears quite
complicated in general, but we can see it directly for the strata of
$\mm_{0,2}({\mathbb P}^r,2)$.  When stratif\/ied according to the dual graphs
of stable maps, $\mm_{0,2}({\mathbb P}^r,2)$ has 9 {\em types} of strata. The
corresponding dual graphs are shown below.

\centerline{\includegraphics{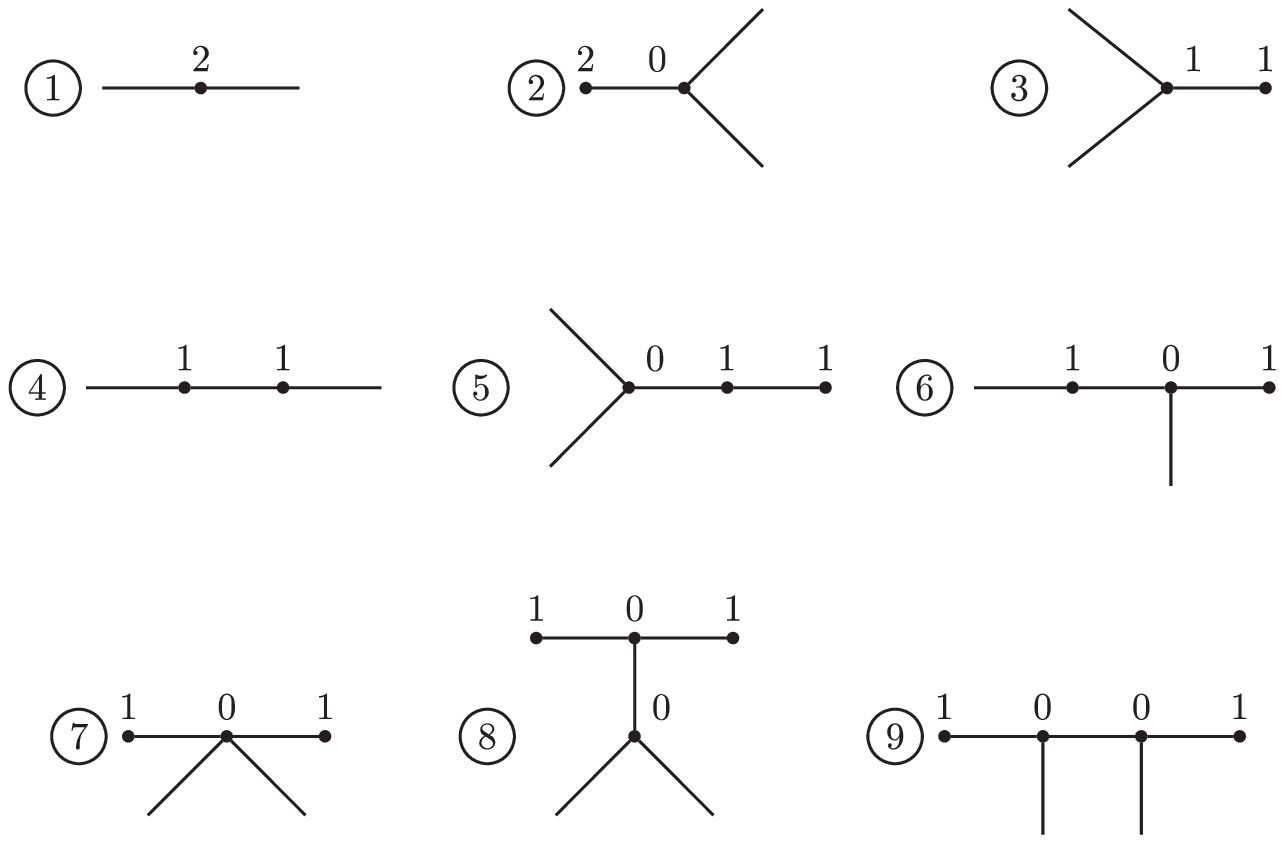}}

All the stratum types are listed, and the assertion clearly
holds in each case.
So we can compute the Serre polynomials of the strata
using Corollary~\ref{serfiber} and Proposition~\ref{Serre00}.

There are actually 10 strata, because there are two strata of type 6 depending
on which marked point is identif\/ied with which tail on the graph.
We use the same numbers to label the strata as
those labeling the corresponding graphs above.
Eight of the strata have no automorphisms, so we can directly compute ordinary
Serre polynomials in these cases. The strata corresponding to
Graphs 7 and 8 have automorphism group $S_2$.
Calculating the $S_2$-equivariant Serre
polynomials of these strata is necessary as
an intermediate step.
We now compute the Serre polynomials of the strata.

Stratum 1 is ${\mathcal{M}}_{0,2}({\mathbb P}^r,2)$. It is an $F({\mathbb P}^1,2)$-bundle over
${\mathcal{M}}_{0,0}({\mathbb P}^r,2)$. Thus Stratum 1 has
Serre polynomial
\[\ser(F({\mathbb P}^1,2))\ser({\mathcal{M}}_{0,0}({\mathbb P}^r,2))
=  (q^2+q)q^{r+1}\frac{[r+1][r]}{[2]} =  q^{r+2}[r+1][r]
.\]

Stratum 2 is isomorphic to the f\/iber product
\[{\mathcal{M}}_{0,1}({{\mathbb P}}^{r},2)\times_{{{\mathbb P}}^{r}}{\mathcal{M}}_{0,3}({{\mathbb P}}^{r},0).\]
Now ${\mathcal{M}}_{0,3}({{\mathbb P}}^{r},0)\simeq{{\mathbb P}}^{r}$, so the Serre polynomial of this stratum
is just
\[\ser({\mathcal{M}}_{0,1}({{\mathbb P}}^{r},2))=\ser({\mathbb P}^1)\ser({\mathcal{M}}_{0,0}({{\mathbb P}}^{r},2))
=q^{r+1}[r+1][r]\]
since ${\mathcal{M}}_{0,1}({{\mathbb P}}^{r},2)$ is a ${\mathbb P}^1$-bundle over ${\mathcal{M}}_{0,0}({{\mathbb P}}^{r},2)$.

Stratum 3 is isomorphic to the f\/iber product
\[{\mathcal{M}}_{0,3}({\mathbb P}^r,1)\times_{{\mathbb P}^r}{\mathcal{M}}_{0,1}({\mathbb P}^r,1)
.\]
The $F({\mathbb P}^1,3)$-bundle ${\mathcal{M}}_{0,3}({\mathbb P}^r,1)$ over ${\mathcal{M}}_{0,0}({\mathbb P}^r,1)$ has
Serre polynomial $(q^3-q)\chews{r+1}{2}$. Similarly,
$\ser({\mathcal{M}}_{0,1}({\mathbb P}^r,1))=(q+1)\chews{r+1}{2}=[r+1][r]$. Thus
Stratum 3
has Serre polynomial
\[\frac{(q^3-q)\chews{r+1}{2}[r+1][r]}{[r+1]}
=(q^2-q)[r+1][r]^2.\]

Stratum 4 is isomorphic to the f\/iber product
\[{\mathcal{M}}_{0,2}({\mathbb P}^r,1)\times_{{\mathbb P}^r}{\mathcal{M}}_{0,2}({\mathbb P}^r,1)
.\]
The $F({\mathbb P}^1,2)$-bundle ${\mathcal{M}}_{0,2}({\mathbb P}^r,1)$  over ${\mathcal{M}}_{0,0}({\mathbb P}^r,1)$
has Serre polynomial $(q^2+q)\chews{r+1}{2}$. Thus Stratum 4
has Serre polynomial
\[\frac{(q^2+q)^2[r+1]^2[r]^2}{[r+1][2]^2}=q^2[r+1][r]^2.\]

Stratum 5 is isomorphic to the f\/iber product
\[{\mathcal{M}}_{0,3}({\mathbb P}^r,0)\times_{{\mathbb P}^r}{\mathcal{M}}_{0,2}({\mathbb P}^r,1)\times_{{\mathbb P}^r}{\mathcal{M}}_{0,1}({\mathbb P}^r,1)
,\]
and this in turn is isomorphic to
${\mathcal{M}}_{0,2}({\mathbb P}^r,1)\times_{{\mathbb P}^r}{\mathcal{M}}_{0,1}({\mathbb P}^r,1)$. So Stratum 5
has Serre polynomial
\[\frac{(q^2+q)\chews{r+1}{2}(q+1)\chews{r+1}{2}}{[r+1]}
=q[r+1][r]^2.\]

A stratum of type 6 is isomorphic to the f\/iber product
\[{\mathcal{M}}_{0,2}({\mathbb P}^r,1)\times_{{\mathbb P}^r}{\mathcal{M}}_{0,3}({\mathbb P}^r,0)\times_{{\mathbb P}^r}{\mathcal{M}}_{0,1}({\mathbb P}^r,1)
.\]
This is isomorphic to  ${\mathcal{M}}_{0,2}({\mathbb P}^r,1)\times_{{\mathbb P}^r}{\mathcal{M}}_{0,1}({\mathbb P}^r,1)$, so
each stratum of type 6 has Serre polynomial
\[\frac{(q^2+q)\chews{r+1}{2}(q+1)\chews{r+1}{2}}{[r+1]}
=q[r+1][r]^2.\]
Thus the total
contribution from strata of type 6 is
\[2q[r+1][r]^2.
\]

Stratum 9 is isomorphic to the f\/iber product
\[{\mathcal{M}}_{0,1}({\mathbb P}^r,1)\times_{{\mathbb P}^r}{\mathcal{M}}_{0,3}({\mathbb P}^r,0)\times_{{\mathbb P}^r}
{\mathcal{M}}_{0,3}({\mathbb P}^r,0)\times_{{\mathbb P}^r}{\mathcal{M}}_{0,1}({\mathbb P}^r,1)
.\]
It
has Serre polynomial
\[\frac{(q+1)^2\chews{r+1}{2}^2}{[r+1]}
=[r+1][r]^2.\]

We now turn our attention to the two strata with automorphisms.
Stratum 8 is isomorphic to the quotient of
\[X={\mathcal{M}}_{0,3}({\mathbb P}^r,0)\times_{{\mathbb P}^r}{\mathcal{M}}_{0,3}({\mathbb P}^r,0)\times_{({\mathbb P}^r)^2}
{\mathcal{M}}_{0,1}({\mathbb P}^r,1)^2\]
by the action of $S_2$.  The f\/irst copy of ${\mathcal{M}}_{0,3}({\mathbb P}^r,0)$
is superf\/luous. The action of
$S_2$ on the cohomology of the second copy of
${\mathcal{M}}_{0,3}({\mathbb P}^r,0)$ is trivial. The action
switches the two factors of ${\mathcal{M}}_{0,1}({\mathbb P}^r,1)$ as well as the two factors
in ${\mathbb P}^r\times{\mathbb P}^r$. Since ${\mathcal{M}}_{0,1}({\mathbb P}^r,1)$ is a f\/iber space
over ${\mathbb P}^r$, we can use Lemma~\ref{ser2X2} and Corollary~\ref{serfiber} in
computing the equivariant Serre polynomial of $X$ to be
\begin{gather*}
 \ser_2({\mathcal{M}}_{0,3}({{\mathbb P}}^{r},0))\ser_2(({\mathcal{M}}_{0,1}({{\mathbb P}}^{r},1)/{{\mathbb P}}^{r})^2)\\
 \qquad{}=  [r+1]\left(\sigma_2\left(\frac{\ser({\mathcal{M}}_{0,1}({{\mathbb P}}^{r},1))}{\ser({{\mathbb P}}^{r})}\right)\one
+\lambda_2\left(\frac{\ser({\mathcal{M}}_{0,1}({{\mathbb P}}^{r},1))}{\ser({{\mathbb P}}^{r})}\right)\epsilon\right) \\
\qquad {}=  [r+1](\sigma_2([r])\one+\lambda_2([r])\epsilon) \\
\qquad {} =  [r+1]\left(\chews{r+1}{2}\one+q\chews{r}{2}\epsilon\right).
\end{gather*}
(As in the proof of Lemma \ref{grasser}, the f\/iber
${\mathcal{M}}_{0,1}({{\mathbb P}}^{r},1)/{{\mathbb P}}^{r}$ is isomorphic to ${\mathbb P}^{r-1}$.)
Now augmentation gives
\[\frac{[r+1]^2[r]}{[2]}\]
as the Serre polynomial of Stratum 8.

Stratum 7 is isomorphic to the quotient of
\[Y={\mathcal{M}}_{0,1}({\mathbb P}^r,1)^2\times_{({\mathbb P}^r)^2}{\mathcal{M}}_{0,4}({\mathbb P}^r,0)\]
by the action of $S_2$, which again switches the squared factors. In addition,
it switches two of the four marked points in ${\mathcal{M}}_{0,4}({\mathbb P}^r,0)$.
Now ${\mathcal{M}}_{0,4}({\mathbb P}^r,0)\simeq M_{0,4}\times{\mathbb P}^r$, and $S_2$ acts trivially
on the ${\mathbb P}^r$ factor. Furthermore,
$M_{0,4}\simeq{\mathbb P}^1\setminus\{0,1,\infty\}$
has Serre polynomial $q-2$. But we need to know $\ser_2(M_{0,4})$
under an $S_2$-action switching
two of the deleted points.
It is not hard to imagine that
$\ser_2(M_{0,4})=(q-1)\one-\epsilon$, but this takes some work to prove.
Considering $M_{0,4}$ as the parameter space of four distinct points
in ${\mathbb P}^1$ modulo automorphisms of ${\mathbb P}^1$, we obtain
$M_{0,4}\simeq F({\mathbb P}^1,4)/\PGL(2)$. Now $\PGL(2)$ acts freely on $F({\mathbb P}^1,4)$.
As a result,
\begin{gather}\label{mm04}
\ser_2(M_{0,4})=\frac{\ser_2(F({\mathbb P}^1,4))}{\ser_2(\PGL(2))}.
\end{gather}
Since the cohomology of $\PGL(2)$ is not af\/fected by the action,
\begin{gather}\label{pgl2}
\ser_2(\PGL(2))=\ser(\PGL(2))=q^3-q.
\end{gather}

We can stratify $({\mathbb P}^1)^4$ into f\/ifteen cells whose closures are
respectively $({\mathbb P}^1)^4$,
the six large diagonals, the seven ``medium diagonals'' where two coordinate
identif\/ications are made, and the small diagonal, so that $F({\mathbb P}^1,4)$
is the complement of the union of all the cells corresponding to diagonals.
We examine how the action
af\/fects cells of each type, subtracting the polynomials for cells that are
removed.
For concreteness, suppose the f\/irst two marked points are switched. Then two
Chow classes in $A^*(({\mathbb P}^1)^4)$ are switched if and only if their dif\/ference is a
multiple of $H_2-H_1$ (where $H_i$ is the standard Chow generator
of $A^*(({\mathbb P}^1)^4)$ obtained by pulling back the hyperplane
class of ${\mathbb P}^1$ under the $i$'th projection), so it is not hard to get
\[
\ser_2(({\mathbb P}^1)^4)=(q^4+3q^3+4q^2+3q+1)\one+(q^3+2q^2+q)\epsilon.
\]
How does the action af\/fect the diagonals removed from $({\mathbb P}^1)^4$?
Exactly two pairs, $(\Delta_{13},\Delta_{23})$ and $(\Delta_{14},\Delta_{24})$, of
the six large diagonals are switched, so the
corresponding cells contribute
\[
(-4q^3+4q)\one+(-2q^3+2q)\epsilon
\]
to the equivariant Serre polynomial, since
these diagonals have been removed.
Exactly two pairs, $(\Delta_{134},\Delta_{234})$ and $(\Delta_{(13)(24)},\Delta_{(14)(23)})$,
among the seven diagonals with two identif\/ications are
switched as well. The corresponding cells contribute
\[
(-5q^2-5q)\one+(-2q^2-2q)\epsilon
\]
to the equivariant Serre polynomial.
The small diagonal is not af\/fected by the action, so it contributes
\[
(-q-1)\one.
\]
Putting these together gives
\[\ser_2(F({\mathbb P}^1,4))=(q^4-q^3-q^2+q)\one+(-q^3+q)\epsilon.\]
Then by (\ref{mm04}) and (\ref{pgl2}), we have the desired result
$\ser_2(M_{0,4})=(q-1)\one-\epsilon$. Using Corollary \ref{serfiber} again,
we thus calculate the equivariant Serre polynomial of $Y$ to be
\begin{gather*}
 \ser_2({\mathcal{M}}_{0,4}({{\mathbb P}}^{r},0))\ser_2(({\mathcal{M}}_{0,1}({{\mathbb P}}^{r},1)/{{\mathbb P}}^{r})^2) \\
\qquad{}= [r+1]((q-1)\one-\epsilon)\left(\chews{r+1}{2}\one+q\chews{r}{2}\epsilon\right) \\
\qquad{}= [r+1]\left(\left((q-1)\chews{r+1}{2}-q\chews{r}{2}\right)\one
+\left((q^2-q)\chews{r}{2}-\chews{r+1}{2}\right)\right)\epsilon.
\end{gather*}
Augmentation gives
\begin{gather*}
[r+1]\left((q-1)\chews{r+1}{2}-q\chews{r}{2}\right)
 =  \frac{[r+1][r]}{[2]}((q-1)[r+1]-q[r-1])\\
\qquad{} =  \frac{[r+1][r]}{[2]}(q^{r+1}+q^r)-\frac{[r+1]^2[r]}{[2]}
 =  [r+1][r]q^r-\frac{[r+1]^2[r]}{[2]}
\end{gather*}
as the Serre polynomial of Stratum 7.

To get the Serre polynomial for the whole moduli space, we add together
the contributions from all the strata
\begin{gather*}
\ser(\mm_{0,2}({\mathbb P}^r,2))\nonumber \\
\quad{}= q^{r+2}[r+1][r]+q^{r+1}[r+1][r]+(q^2-q)[r+1][r]^2
+q^2[r+1][r]^2+q[r+1][r]^2\nonumber \\
\quad\phantom{=}{} +2q[r+1][r]^2+[r+1][r]^2+\frac{[r+1]^2[r]}{[2]}
+[r+1][r]q^r-\frac{[r+1]^2[r]}{[2]}\nonumber \\
\quad{}= [r+1][r](q^{r+2}+q^{r+1}+(q^2-q)[r]+q^2[r]+3q[r]+[r]+q^r)\nonumber \\
\quad{}= [r+1][r](q^{r+2}+q^{r+1}+q^r+[r](2q^2+2q+1))\nonumber \\
\quad{}= [r+1][r]\left(q^{r+2}+q^{r+1}+q^r+2\sum_{i=2}^{r+1}q^i+2\sum_{i=1}^{r}q^i
+\sum_{i=0}^{r-1}q^i\right)\nonumber \\
\quad{}= \left(\sum_{i=0}^r q^i\right)\left(\sum_{i=0}^{r-1} q^i\right)
\left(\sum_{i=0}^{r+2} q^i+2\sum_{i=1}^{r+1} q^i+2\sum_{i=2}^{r} q^i\right)
.
\end{gather*}
Evaluating this sum at $q=1$ gives the Euler characteristic
$(r+1)r(5r+3)$.
\end{proof}

\subsection[Formulas for the Betti numbers of $\mm_{0,2}({\mathbb P}^r,2)$]{Formulas for the Betti numbers of $\boldsymbol{\mm_{0,2}({\mathbb P}^r,2)}$}
\label{sec:betform}

Let $\alpha_i$ denote the $i$'th Betti number of the f\/lag variety
${\mathbb F}(0,1;r)$ of point-line pairs in ${\mathbb P}^r$ such that the point lies
on the line. Recall from the proof of Lemma \ref{grasser} that
$\ser({\mathbb F}(0,1;r))=[r+1][r]$. The product $[r+1][r]$ also appears as
a factor in the Serre polynomial (\ref{ser2r2}) of
$\mm_{0,2}({\mathbb P}^r,2)$, making its coef\/f\/icients especially relevant to
our computations. It is easy to see that the Betti numbers  of
${\mathbb F}(0,1;r)$ initially follow the pattern $(1, 2, 3, \dots)$,
so that for the f\/irst half of the Betti numbers we have $\alpha_i=i+1$.
Since $\dim{\mathbb F}(0,1;r)=2r-1$ is always odd, it always has an even number
of Betti numbers.  By Poincar\'{e} duality, it
follows that the middle two Betti numbers are both~$r$, and the Betti numbers
then decrease back to~1.
It can be checked that all the Betti numbers are given by the
formula{\samepage
\[\alpha_i=r+\frac{1}{2}-\left|r-\frac{1}{2}-i\right|\]
for $i\in{0}\cup\underline{2r-1}$, and $\alpha_i=0$ otherwise.}

Let $\beta_j$ be the $j$'th Betti number of $\mm_{0,2}({\mathbb P}^r,2)$.  By distributing over the rightmost
set of parentheses in Equation \ref{ser2r2}, we can reduce the computation of $\beta_j$ to f\/inding coef\/f\/icients
of expressions of the form $[r+1][r][m]$, where $m\in\{r-1,r+1,r+3\}$. But these can be expressed in terms of the
$\alpha_i$, and in this way we get the following formulas for the Betti numbers of
$\mm_{0,2}({\mathbb P}^r,2)$:
\begin{gather*}
\beta_{j}=
\begin{cases}
\displaystyle \sum_{i=0}^{j}\alpha_i+2\sum_{i=0}^{j-1}\alpha_i+2\sum_{i=0}^{j-2}\alpha_i
&\text{ if } j\leq r, \vspace{1mm}\\
\displaystyle  \sum_{i=0}^{r+1}\alpha_i+2\sum_{i=0}^{r}\alpha_i+2\sum_{i=1}^{r-1}\alpha_i
&\text{ if } j=r+1,\vspace{1mm}\\
\displaystyle  \sum_{i=j-r-2}^{j}\alpha_i+2\sum_{i=j-r-1}^{j-1}\alpha_i+2\sum_{i=j-r}^{j-2}
\alpha_i
&\text{ if } r+2\leq j\leq 2r-1, \vspace{1mm}\\
\displaystyle  \sum_{i=r-2}^{2r-1}\alpha_i+2\sum_{i=r-1}^{2r-1}\alpha_i+2\sum_{i=r}^{2r-2}
\alpha_i   &\text{ if } j=2r,\vspace{1mm}\\
\displaystyle \sum_{i=j-r-2}^{2r-1}\alpha_i+2\sum_{i=j-r-1}^{2r-1}\alpha_i
+2\sum_{i=j-r}^{2r-1}\alpha_i   &\text{ if } 2r+1\leq j\leq 3r+1.
\end{cases}
\end{gather*}

We can come up with an especially explicit description of $\beta_j$ for
$j<r$ since we know $\alpha_i=i+1$ for $i<r$. Also, $\alpha_r=r$, which gives the
second part below.

\begin{corollary}
\begin{enumerate}\itemsep=0pt
\item For $j<r$, the $j$'th Betti number of $\mm_{0,2}({\mathbb P}^r,2)$ is
\[\beta_j=\frac{5}{2}j^2+\frac{3}{2}j+1.\]

\item Furthermore
\[\beta_r=\frac{5}{2}r^2+\frac{3}{2}r.\]

\end{enumerate}
\end{corollary}

As a consequence of this, a particular Betti number of $\mm_{0,2}({\mathbb P}^r,2)$
stabilizes as $r$ becomes large.

\begin{corollary}
For all $r>j$, the $j$'th Betti number of $\mm_{0,2}({\mathbb P}^r,2)$ is
$\beta_j=\frac{5}{2}j^2+\frac{3}{2}j+1$.
\end{corollary}

Let $\bar{\beta}_j$ be this limiting value. We have
\[\bar{\beta}_0=1,\quad
\bar{\beta}_1=5,\quad
\bar{\beta}_2=14,\quad
\bar{\beta}_3=28,\quad
\bar{\beta}_4=47,\quad
\bar{\beta}_5=71,\quad\ldots.\]

\subsection[Poincar\'{e} polynomials of $\mm_{0,1}({\mathbb P}^r,2)$ and
$\mm_{0,2}({\mathbb P}^r,2)$
for small $r$]{Poincar\'{e} polynomials of $\boldsymbol{\mm_{0,1}({\mathbb P}^r,2)}$ and
$\boldsymbol{\mm_{0,2}({\mathbb P}^r,2)}$
for small $\boldsymbol{r}$}\label{sec:poincaresmallr}

Using the same procedure as above,
one can easily compute the Poincar\'{e} polynomial of
$\!\mm_{0,1}\!({\mathbb P}^r\!,2),\!\!$ which is also needed in the sequel~\cite{C}.

\begin{proposition}
If $r$ is even, the Poincar\'{e} polynomial of $\mm_{0,1}({\mathbb P}^r,2)$ is
\[\ser(\mm_{0,1}({\mathbb P}^r,2))=
\left(\sum_{i=0}^r q^i\right)\left(\sum_{i=0}^{(r-2)/2} q^{2i}\right)
\left(\sum_{i=0}^{r+2} q^i+\sum_{i=1}^{r+1} q^i+\sum_{i=2}^{r} q^i\right)
,\]
and if $r$ is odd, the Poincar\'{e} polynomial of $\mm_{0,1}({\mathbb P}^r,2)$ is
\[\ser(\mm_{0,1}({\mathbb P}^r,2))=
\left(\sum_{i=0}^{r-1} q^i\right)\left(\sum_{i=0}^{(r-1)/2} q^{2i}\right)
\left(\sum_{i=0}^{r+2} q^i+\sum_{i=1}^{r+1} q^i+\sum_{i=2}^{r} q^i\right)
.\]
\end{proposition}

Thus, for small values of $r$, we get the explicit Poincar\'{e}
polynomials listed in Tables~\ref{ep01r2} and~\ref{ep02r2}.

\begin{table}[h]\centering
\caption{Euler characteristics and Poincar\'{e} polynomials for
$X=\mm_{0,1}({\mathbb P}^r,2)$. \label{ep01r2}}

\vspace{0.5mm}

\begin{tabular}{|rrl|}
\hline
$r$ & $\chi(X)$ & $\ser(X)$ \\
1 &  6 &  $1+2q+2q^2+q^3$ \\
2 & 27 & $1+3q+6q^2+7q^3+6q^4+3q^5+q^6$ \\
3 & 72 & $1+3q+7q^2+11q^3+14q^4+14q^5+11q^6+7q^7+3q^8+q^9$ \\
4 & 150 & $1+3q+7q^2+12q^3+18q^4+22q^5+24q^6$ \\
  & & \hspace{0.1in} $+22q^7+18q^8+12q^9+7q^{10}+3q^{11}+q^{12}$ \\
5 & 270 & $1+3q+7q^2+12q^3+19q^4+26q^5+32q^6+35q^7$ \\
& & \hspace{0.1in} $+35q^8+32q^9+26q^{10}+19q^{11}+12q^{12}+7q^{13}+3q^{14}+q^{15}$ \\
6 & 441 & $1+3q+7q^2+12q^3+19q^4+27q^5+36q^6+43q^7+48q^8+49q^9$ \\
& & \hspace{0.1in} $+48q^{10}+43q^{11}+36q^{12}+27q^{13}+19q^{14}+12q^{15}+7q^{16}+3q^{17}+q^{18}$ \\
7 & 672 & $1+3q+7q^2+12q^3+19q^4+27q^5+37q^6+47q^7$ \\
& & \hspace{0.1in}$+56q^8+62q^9+65q^{10}+65q^{11}+62q^{12}+56q^{13}+47q^{14}$ \\
& & \hspace{0.1in} $+37q^{15}+27q^{16}+19q^{17}+12q^{18}+7q^{19}+3q^{20}+q^{21}$ \\
\hline
\end{tabular}\vspace{-3mm}
\end{table}

\begin{table}[h]\centering
\caption{Euler characteristics and Poincar\'{e} polynomials for
$Y=\mm_{0,2}({\mathbb P}^r,2)$. \label{ep02r2}}

\vspace{0.5mm}

\begin{tabular}{|rrl|}
\hline
$r$ & $\chi(Y)$ & $\ser(Y)$ \\
1 & 16 & $1+4q+6q^2+4q^3+q^4$  \\
2 & 78 & $1+5q+13q^2+20q^3+20q^4+13q^5+5q^6+q^7$   \\
3 & 216 & $1+5q+14q^2+27q^3+39q^4+44q^5+39q^6+27q^7+14q^8+5q^9+q^{10}$  \\
4 & 460 & $1+5q+14q^2+28q^3+46q^4+63q^5+73q^6+73q^7$ \\
   & & \hspace{0.1in} $+63q^8+46q^9+28q^{10}+14q^{11}+5q^{12}+q^{13}$   \\
5 & 840 & $1+5q+14q^2+28q^3+47q^4+70q^5+92q^6+107q^7+112q^8$ \\
& & \hspace{0.1in} $+107q^9+92q^{10}+70q^{11}+47q^{12}+28q^{13}+14q^{14}+5q^{15}+q^{16}$   \\
6 & 1386 & $1+5q+14q^2+28q^3+47q^4+71q^5+99q^6+126q^7$ \\
& & \hspace{0.1in} $+146q^8+156q^9+156q^{10}+146q^{11}+126q^{12}+99q^{13}$ \\
& & \hspace{0.1in} $+71q^{14}+47q^{15}+28q^{16}+14q^{17}+5q^{18}+q^{19}$   \\
7 & 2128 & $1+5q+14q^2+28q^3+47q^4+71q^5+100q^6+133q^7+165q^8$ \\
& & \hspace{0.1in} $+190q^9+205q^{10}+210q^{11}+205q^{12}+190q^{13}+165q^{14}+133q^{15}$ \\
& & \hspace{0.1in} $+100q^{16}+71q^{17}+47q^{18}+28q^{19}+14q^{20}+5q^{21}+q^{22}$   \\
\hline
\end{tabular}\vspace{-2mm}
\end{table}

\section[An additive basis for $\mm_{0,2}({\mathbb P}^r,2)$]{An additive basis for $\boldsymbol{\mm_{0,2}({\mathbb P}^r,2)}$}
\label{sec:basis}

A presentation for the Chow rings $A^*({\mathcal{M}}_{0,0}({\mathbb P}^r,d))$ is described in \cite{P4} (see also \cite{BO}).
Naturally then, an additive basis
for these rings is readily available. In this section, we describe an additive basis for $A^*(\mm_{0,2}({\mathbb P}^r,2))$
in terms of the additive bases for $A^*({\mathcal{M}}_{0,0}({\mathbb P}^r,1))$ and $A^*({\mathcal{M}}_{0,0}({\mathbb P}^r,2))$ using a
decomposition of $\mm_{0,2}({\mathbb P}^r,2)$ and excision. (Note that an additive basis has
now been more explicitly described in \cite{MM2}.)

Let $X$ be the locus in $\mm_{0,2}({\mathbb P}^r,2)$ where the curve has two degree one components (with degree zero
components also allowed). It is a divisor, and
its complement is the open locus $U$ where the curve has a degree two component. In the notation of Section
\ref{sec:ser},
$U$ is the union of Strata 1 and 2.

Note that $U$ is a ${\mathbb P}^1\times{\mathbb P}^1$-bundle over ${\mathcal{M}}_{0,0}({\mathbb P}^r,2)$. Let $H_1$ and $H_2$ be the two hyperplane divisor
classes in ${\mathbb P}^1\times{\mathbb P}^1$. Then any class in $A^k(U)$ can be expressed in the form $\alpha+\beta_1 H_1+\beta_2 H_2+ \gamma H_1H_2$, where
$\alpha\in A^k({\mathcal{M}}_{0,0}({\mathbb P}^r,2))$, $\beta_i\in A^{k-1}({\mathcal{M}}_{0,0}({\mathbb P}^r,2))$, and
$\gamma\in A^{k-2}({\mathcal{M}}_{0,0}({\mathbb P}^r,2))$. Thus we can write
\begin{gather*}
A^k(U)\simeq
A^k({\mathcal{M}}_{0,0}({\mathbb P}^r,2))\oplus A^{k-1}({\mathcal{M}}_{0,0}({\mathbb P}^r,2))\oplus A^{k-1}({\mathcal{M}}_{0,0}({\mathbb P}^r,2))
\oplus A^{k-2}({\mathcal{M}}_{0,0}({\mathbb P}^r,2)).
\end{gather*}

Consider the exact sequence
\[A^{k-1}(X)\longrightarrow A^k(\mm_{0,2}({\mathbb P}^r,2))\longrightarrow A^k(U)\longrightarrow 0
.\]
This says $A^k(\mm_{0,2}({\mathbb P}^r,2))$ is the direct sum of $A^k(U)$ and the image of $A^{k-1}(X)$.
We have reduced to studying the latter.
Let $Y$ be the sublocus of $X$
where there are no marked points on the degree one components. Then $Y$ is a divisor in $X$, and
we can consider the exact sequence
\[A^{k-2}(Y)\longrightarrow A^{k-1}(X)\longrightarrow A^{k-1}(V)\longrightarrow 0
,\]
where $V=X-Y$. We can further decompose $V$ into $V_1\coprod V_2$, where $V_2$ contains the locus where each degree one component has a marked point. In $V_2$ we also allow the second marked point to approach the node. In the notation of
Section~\ref{sec:ser}, $V_2$ is the union of Stratum 4 and one of the strata of type~6. Thus $V_2$ is a ${\mathbb P}^1\times{\mathbb A}^1$-bundle over the boundary divisor
$D\simeq\mm_{0,1}({\mathbb P}^r,1)\times_{{\mathbb P}^r}\mm_{0,1}({\mathbb P}^r,1)$ in $\mm_{0,0}({\mathbb P}^r,2)$.
Since the f\/iber product is a ${\mathbb P}^{r-1}$-bundle over $\mm_{0,1}({\mathbb P}^r,1)$,
and $\mm_{0,1}({\mathbb P}^r,1)$ is a ${\mathbb P}^1$-bundle
over $\mm_{0,0}({\mathbb P}^r,1)\simeq{\mathcal{M}}_{0,0}({\mathbb P}^r,1)$, we f\/ind overall that $D$ is a
${\mathbb P}^{r-1}\times{\mathbb P}^1$-bundle over ${\mathcal{M}}_{0,0}({\mathbb P}^r,1)$.
Similarly,
$V_1$, which is the union of Stratum 3, Stratum 5, and the other stratum of type 6, is also
a ${\mathbb P}^1\times{\mathbb A}^1$-bundle over $D$, and thus a ${\mathbb P}^{r-1}\times({\mathbb P}^1)^2\times{\mathbb A}^1$-bundle
over ${\mathcal{M}}_{0,0}({\mathbb P}^r,1)$. We have
\begin{gather*}
A^{k-1}(V_2)\simeq A^{k-1}(V_1)
 \simeq  \oplus_{i=0}^{r-1}\left(A^{k-i-1}({\mathcal{M}}_{0,0}({\mathbb P}^r,1))\oplus A^{k-i-2}({\mathcal{M}}_{0,0}({\mathbb P}^r,1))\right. \\
\left.\phantom{A^{k-1}(V_2)\simeq A^{k-1}(V_1)
 \simeq}{} \oplus A^{k-i-2}({\mathcal{M}}_{0,0}({\mathbb P}^r,1))\oplus A^{k-i-3}({\mathcal{M}}_{0,0}({\mathbb P}^r,1))\right).
\end{gather*}

Finally, $Y$ is the union of Strata 7, 8, and 9, and is a ${\mathbb P}^1$-bundle over the codimension two boundary stratum $Z$ in $\mm_{0,1}({\mathbb P}^r,2)$. (This is the
boundary locus where the domain curves have three components.) Now $Z$ is isomorphic to the $S_2$-quotient of
$\mm_{0,1}({\mathbb P}^r,1)^2\times_{({\mathbb P}^r)^2}\mm_{0,3}({\mathbb P}^r,0)$ that arises by switching the factors in the squares. This f\/iber
product is isomorphic to a $({\mathbb P}^{r-1})^2$-bundle over ${\mathbb P}^r$. So
\[A^{k-2}(Y)\simeq\left(\oplus_{i=0}^{r-1}\oplus_{j=0}^{r-1}A^{k-i-j-2}({\mathbb P}^r)\right)^{S_2}
\oplus \left(\oplus_{i=0}^{r-1}\oplus_{j=0}^{r-1}A^{k-i-j-3}({\mathbb P}^r)\right)^{S_2}.
\]

Putting all this together, we attain an additive basis for $A^*(\mm_{0,2}({\mathbb P}^r,2))$. It is given in degree~$k$ by
\begin{gather*}
 A^k(\mm_{0,2}({\mathbb P}^r,2))  \simeq
 A^k({\mathcal{M}}_{0,0}({\mathbb P}^r,2))\oplus A^{k-1}({\mathcal{M}}_{0,0}({\mathbb P}^r,2))\nonumber \\
\phantom{A^k(\mm_{0,2}({\mathbb P}^r,2))  \simeq}{} \oplus A^{k-1}({\mathcal{M}}_{0,0}({\mathbb P}^r,2))\oplus A^{k-2}({\mathcal{M}}_{0,0}({\mathbb P}^r,2))\nonumber \\
\phantom{A^k(\mm_{0,2}({\mathbb P}^r,2))  \simeq}{} \oplus \oplus_{i=0}^{r-1}\left(A^{k-i-1}({\mathcal{M}}_{0,0}({\mathbb P}^r,1))\oplus A^{k-i-2}({\mathcal{M}}_{0,0}({\mathbb P}^r,1))\right.\nonumber  \\
 \left. \phantom{A^k(\mm_{0,2}({\mathbb P}^r,2))  \simeq}{}\oplus A^{k-i-2}({\mathcal{M}}_{0,0}({\mathbb P}^r,1))\oplus A^{k-i-3}({\mathcal{M}}_{0,0}({\mathbb P}^r,1))\right)\nonumber  \\
\phantom{A^k(\mm_{0,2}({\mathbb P}^r,2))  \simeq}{} \oplus \oplus_{i=0}^{r-1}\left(A^{k-i-1}({\mathcal{M}}_{0,0}({\mathbb P}^r,1))\oplus A^{k-i-2}({\mathcal{M}}_{0,0}({\mathbb P}^r,1))\right.\nonumber  \\
\phantom{A^k(\mm_{0,2}({\mathbb P}^r,2))  \simeq}{} \oplus A^{k-i-2}\left.(
{\mathcal{M}}_{0,0}({\mathbb P}^r,1))\oplus A^{k-i-3}({\mathcal{M}}_{0,0}({\mathbb P}^r,1))\right)\nonumber  \\
\phantom{A^k(\mm_{0,2}({\mathbb P}^r,2))  \simeq}{} \oplus \left(\oplus_{i=0}^{r-1}\oplus_{j=0}^{r-1}A^{k-i-j-2}({\mathbb P}^r)\right)^{S_2}
\oplus \left(\oplus_{i=0}^{r-1}\oplus_{j=0}^{r-1}A^{k-i-j-3}({\mathbb P}^r)\right)^{S_2}.
\end{gather*}

Of course, this is only a generating set {\em a priori}, but by comparing with the Betti numbers we can see that the
generators are independent. In other words, we need only
verify that the expression above
gives the same Serre polynomial for $A^*(\mm_{0,2}({\mathbb P}^r,2))$ as found in Section
\ref{sec:ser}. Keeping in
mind that the Serre polynomial grades by dimension rather than codimension, and using the notation of Section
\ref{sec:ser}, from the decomposition above we obtain
\begin{gather*}
 (q+1)^2q^{r+1}\chews{r+1}{2}+2\sum_{i=0}^{r-1}q^{i+1}(q+1)^2\chews{r+1}{2}+(q+1)[r+1]\sigma_2([r]) \\
\qquad{}= [r+1][r]((q+1)q^{r+1}+2q[r](q+1)+[r+1])\\
\qquad{}= [r+1][r]\left(q^{r+2}+q^{r+1}+2\sum_{i=1}^rq^i+2\sum_{i=2}^{r+1}q^i+\sum_{i=0}^rq^i\right) \\
\qquad{}= [r+1][r]\left(\sum_{i=0}^{r+2}q^i+2\sum_{i=1}^{r+1}q^i+2\sum_{i=2}^rq^i\right),
\end{gather*}
in agreement with equation~\eqref{ser2r2}.

\subsection*{Acknowledgements}

The content of this article derives from a part of my doctoral dissertation at Oklahoma State University.
I am deeply grateful to my dissertation adviser, Sheldon Katz for f\/inancial support, insight,
encouragement, and inspiration.
William Jaco and
Alan Adolphson provided
additional funding
during work on this project.
I appreciate the hospitality of the University of Illinois
mathematics department during my years
as a visiting graduate student there. I also acknowledge
with gratitude the
Oklahoma State University mathematics department for extended
support during that time.

\pdfbookmark[1]{References}{ref}
\LastPageEnding

\end{document}